\documentclass[11pt]{amsart}
\font\bbbld=msbm10 scaled\magstephalf

\newcommand{\bfR}{\hbox{\bbbld R}}

\newcommand{\vp}{{\bf p}}
\newcommand{\ol}{\overline}

\newtheorem{theorem}{Theorem}[section]
\newtheorem{lemma}[theorem]{Lemma}

\newtheorem{corollary}[theorem]{Corollary}
\newtheorem{notation}[theorem]{Notation}

\theoremstyle{definition}

\theoremstyle{remark}
\newtheorem{remark}[theorem]{Remark}

\numberwithin{equation}{section}

\begin{document}
\setlength{\baselineskip}{1.2\baselineskip}

\centerline{\bf Hessian equations with Infinite Dirichlet Boundary
Value}
 \vskip 1cm
 \centerline{Huaiyu Jian    }
 \centerline{Department of  Mathematical Sciences, Tsinghua University,
Beijing 100084, China}

\vskip 2cm

\noindent{\bf Correspondence Address }

Huaiyu Jian

  Department of  Mathematical Sciences

 Tsinghua University

  Beijing 100084, P. R. China

\noindent {\bf Email}\ \ hjian@math.tsinghua.edu.cn

\noindent {\bf Telfax}  \ \ 86-10-62781785

\noindent {\bf Tel}\ \ 86-10-62772864

 \newpage

\title[Hessian equations]
{Hessian equations with Infinite Dirichlet Boundary Value}
\author{Huaiyu Jian}
\address{Department of  Mathematical Sciences, Tsinghua University,
Beijing 100084, China} \email{hjian@math.tsinghua.edu.cn}

\thanks{ The research was supported by
National 973-Project from MOST and   Basic Research Grant of Tsinghua University. \\
The work was  partially done while the   author was visiting the
Institute for Mathematical Sciences, National University of
Singapore in 2004. The visit was supported by the Institute.
 }


\maketitle



\noindent {\bf Abstract}   In this paper, we study the Hessian
equation with infinite Dirichlet (blow-up) boundary value
conditions. Using radial functions and  techniques of ordinary
differential inequality, we construct various barrier functions
(super-solution and sub-solution). Existence and non-existence
theorems are proved by those barriers, maximum principle and
theory of viscous solutions. Furthermore, generic boundary blow-up
rates for the solutions are derived.



\vskip 0.5cm
\section{Introduction}
\setcounter{equation}{0}

  Let $n$ and $1\leq k\leq n$ be positive integers. Recall that the
$k$-th elementary symmetric function for $\lambda\in\bfR^n$ is
defined by
$$ S_{k}(\lambda) =S_{k}(\lambda_{1},...,\lambda_{n})
=\sum_{{1\leq{i_{1}}<{i_{2}}<...<{i_{k}}\leq{n}}}\lambda_{i_{1}}\lambda_{i_{2}}...\lambda_{i_{k}}$$
and the $k$-th elementary symmetric function over the space $
M_{s}(\bfR^{n})$ of all the $ n\times n $ real symmetric matrices
is given by
$$ S_{k}(\mathbf{A})=S_{k}(\lambda_{1},...,\lambda_{n}), \;\;\forall \;\mathbf{A}\in M_{s}(\bf R^n) ,$$
where $ (\lambda_{1},...,\lambda_{n}) $ are the eigenvalues of the
matrix $\mathbf{A}$.

 Let $\Omega$ be a domain in $\bfR^{n}$
and $\psi$ a positive function defined on $\Omega \times \bfR
\times \bfR^{n}$. In this paper, we deal with the Hessian equation
\begin{equation}
S_{k}(D^2 u)  =  \psi (x, u, Du) \;\; \mbox{in $\Omega$}
\end{equation}
with the singular boundary value condition
\begin{equation}
 u(x) = + \infty \;\; \mbox{on $\partial \Omega$}.
\end{equation}
Here $ D^2 u=(u_{ij})_{n\times n}$ is the Hessian matrix of $u$
and (1.2) means that $u(x) \rightarrow+ \infty$ as
$d(x)=d(x,\partial \Omega)\rightarrow0$, where $d(x,\partial
\Omega)$ denotes the distance of  point $x\in\Omega$ from
$\partial \Omega$.

Obvious examples of Hessian in (1.1) are Laplace operator, $k=1$,
and the Monge-Amp\`ere operator, $k=n$. General Hessian operators
were studied by many authors. See, for example, [2, 5, 6, 16, 18]
for Hessian in $\Omega$ and [1, 8, 9, 10] for Hessian on the
sphere.

A natural class of functions for the solutions to (1.1)-(1.2) is
k-convex functions. Recall that a function $u\in C^{2}(\Omega)$ is
called k-convex (or strictly k-convex) if
$(\lambda_{1},...,\lambda_{n})\in\overline{\Gamma_{k}}(or(\lambda_{1},...,\lambda_{n})\in\Gamma_{k})$
for every $x\in\Omega$, where $\lambda_{1},...,\lambda_{n}$ are
the eigenvalue of $ D^2 u(x)$ and $\Gamma_{k}$ is the connected
component \{${\lambda\in \bfR^{n}:S_{k}(\lambda)>0}$\} containing
the positive cone
$$\Gamma^{+}=\{\lambda=(\lambda_{1},...,\lambda_{n})\in
\bfR^{n}:\lambda_{i}>0,i=1,2,...,n\}.$$ It follows from [1] that
\begin{equation}
\Gamma^{+}=\Gamma_{n}\subset...\subset\Gamma_{k+1}\subset\Gamma_{k}\subset...\subset\Gamma_{1}
\end{equation}
and $S_{k}(D^2 u)$ turns to be elliptic in the class of k-convex
functions.

The problem (1.1)-(1.2) was studied in [11, 14] for $k=1$, Laplace
operator, and in [3, 4, 12, 13] for $k=n$, Monge-Amp\`ere
operator, and $\psi=\psi(x,u)$ independent of $Du$. The results of
Matero [13] were extended by Salani [15] to some Hessian equation;
while the results of Cheng and Yau [3, 4] were generalized
recently by Guan and Jian [7], in which various existence and
non-existence results were shown for rather general
$\psi=\psi(x,u,Du)$ and the optimal growth condition of
$\psi(x,z,\vp)$ was given for the existence of (1.1)-(1.2) in the
case $k=n$. The aim of this paper is to extend the main results of
[7] to the case $k\in\{1,2,...,n-1\}$. The difficulty here arises
when one tries to construct barriers which is necessary for the
existence or non-existence of problem (1.1)-(1.2). The methods and
results of this paper are different from those of [15].

From now on,we assume
\begin{equation}
k\in\{1,2,...,n-1\}.
\end{equation}

Our main results are stated as follows.
\begin{theorem}
 Let $\Omega$ be a bounded domain in $\bfR^{n}$. If there are
 constants $M,\gamma,q \geq 0,\gamma+q \leq k$, such that
\begin{equation}
0 \leq \psi(x,z,\vp) \leq M(1+(z^{+})^{q})(1+{|\vp|}^{\gamma}),
\;\;\forall \;(x,z,\vp)\in \Omega\times\bfR\times\bfR^{n},
\end{equation}
where $z^{+}=max\{z,0\}$, then there exists no k-convex solution
to (1.1)-(1.2) in $C^{2}(\Omega)$.
\end{theorem}
\begin{theorem}
If there are constants $\alpha>1$ and $M>0$ such that
\begin{equation}
 \psi(x,z,\vp)\geq
M(1+{|\vp|}^{k})^{\alpha},\;\;\forall \;(x,z,\vp)\in
\Omega\times\bfR\times\bfR^n,
\end{equation}
and $\Omega$ is a domain containing some ball of radius a where
$$ a>\big[\frac{(k+1)(n-1)!}{M(\alpha-1)k!(n-k-1)!}\big]^{\frac{1}{k}},$$
then there exists no k-convex solution to (1.1)-(1.2) in
$C^{2}(\Omega)$.
\end{theorem}

We will deal with the existence of problem (1.1)-(1.2) in
viscosity sense. For the details of viscosity solutions to Hessian
equations like (1.1), as for the notion of k-convexity in
viscosity sense, we refer to [18].
\begin{theorem}
Let $\Omega$ be a bounded strictly convex domain in $\bfR^{n}$ and
$\psi \in C^{\infty} (\ol{\Omega} \times \bfR \times \bfR^{n})$
satisfy
\begin{equation}
\psi(x,z,\vp)>0,\;\psi_{z}(x,z,\mathbf{p})>0,\;\;\forall
\;(x,z,\vp)\in \Omega\times\bfR\times\bfR^{n}.
\end{equation}
Suppose that there exist $q>k$ and $M>0$ such that for all $
(x,z,\vp)\in \Omega\times\bfR\times\bfR^{n},$
\begin{equation} \psi(x,z,\vp)\geq
M(z^{+})^{q}
\end{equation}
and
\begin{equation}
\psi(x,z,\vp)\leq\phi(z)(1+|\vp|^{k}),
\end{equation}
where $\phi\in C^1(R^n)$ is a  positive nondecreasing function
satisfying
\begin{equation}
\sup_{z\leq 0}e^{-\varepsilon z}\phi(z)<+\infty
\end{equation}
for some $\varepsilon>0$, and for each fixed
$(x,z)\in\overline{\Omega}\times\bfR,$
\begin{equation}
\psi^{\frac{1}{k}}(x,z,\vp) \;\mbox{is convex in $\vp$},\;
\inf_{\vp \in R^n}\psi>0,\; \sup_{\vp \in R^n
}\frac{|D_x(\psi^{\frac{1}{k}})|}{1+|\vp|}<+\infty.
\end{equation}
Then there exists a k-convex viscosity solution $u\in
C^{0,1}_{loc}(\Omega)$ to (1.1)-(1.2). Moreover, there exist
functions $\underline{h},\;\overline{h}\in C({\bfR}^{+})$ with
$\underline{h}(r),\;\overline{h}(r)\rightarrow +\infty$ as
$r\rightarrow 0,$ such that
\begin{equation}
\underline{h}(d(x))\leq u(x)\leq\overline{h}(d(x)), \forall x\in
\Omega.
\end{equation}
\end{theorem}

\begin{remark}
One can see easily that there is a large kind of functions
satisfying (1.7)-(1.11). For example,
$$\psi (x,
z,\vp)=(a(x)+b(x)e^{rz})(1+\sum_{i=1}^{n}p_i^{2k})^{\frac{1}{2}}$$
or
$$\psi (x,
z,\vp)=(a(x)+b(x)\eta
(z))(1+\sum_{i=1}^{n}p_i^{2k})^{\frac{1}{2}},$$ where $a,b$ are
positive smooth functions on $\ol{\Omega},$ $r$ is a positive
constant, and $\eta$ is a nonnegative strictly increasing smooth
function satisfying $\eta (z)\geq |z|^q $ for all $z\geq 1$ and
for some constant $q>k.$
\end{remark}

\begin{remark}
It seems possible to get better regularity than $C^{0,1}$ for the
solution obtained in theorem 1.3.  One helpful approach is to use
the methods on gradient estimates in [5].
\end{remark}

\section{A comparison principle and uniqueness}
\setcounter{equation}{0}

This section is similar to section 2 in [7]. For the sake of
convenience, we will give the details here.

Suppose that $u,v\in C^{2}(\Omega)$ are $k$-convex satisfying
\begin{equation}
S_{k}(D^2 u)  \geq  \psi (x, u, Du) \; \mbox{and $S_{k}(D^2 v)\leq
\phi(x, v, Dv)$} \;\;\mbox{in $\Omega$,}
\end{equation}
where $\phi,\psi\in C^{1} (\Omega \times \bfR \times \bfR^{n})$,
such that
\begin{equation}
\psi (x, z,\vp)\geq \phi(x, z,\vp), \;\;\forall \;(x,
z,\vp)\in(\Omega \times \bfR \times \bfR^{n}).
\end{equation}
We will see if $u\leq v$ in $\Omega$. Let $M^{k}_{s}(\bfR^{n})$ be
the subset of $M_{s}(\bfR^{n})$ , with the eigenvalues
$\lambda_{1},...,\lambda_{n}$  satisfying
$(\lambda_{1},...,\lambda_{n})\in \Gamma_{k} $. Recall that
\begin{equation}(\frac{{\partial S}_{k}^{\frac {1}{k}}(\mathbf{w})}{{\partial
  {\mathbf{w}}}_{ij}})_{n \times n}>0 \; (\mbox{or $\geq 0$}), \;\; \forall
 \;\mathbf{w} \in M^{k}_{s}(\bfR^n) \; (\mbox{or
 $\overline{M^{k}_{s}(\bfR^{n})}$})
\end{equation}
(See [2, 18]). Hence for any $ \mathbf{w}_{1},\mathbf{w}_{2} \in
M^{k}_{s}(\bfR^{n}) $ (or $ \overline{M^{k}_{s}(\bfR^{n})}$),
\begin{equation}
S_{k}(\mathbf{w}_{1})> S_{k}(\mathbf{w}_{2}) \;(\mbox{or $
S_{k}(\mathbf{w}_{1})\geq S_{k}(\mathbf{w}_{2})$})
\end{equation}
if $\mathbf{w}_{1}-\mathbf{w}_{2}$ is positive definite (or
semi-definite).
\begin{lemma}
Suppose that $\Omega\subset\bfR^{n}$ is a bounded domain, and
$u,v\in C(\overline{\Omega})$ satisfy $u\leq v$ on $\partial
\Omega$. If either $\psi_{z}(x,z,\vp)>0$ or $\phi_{z}(x,z,\vp)>0$
for any $(x, z,\vp)\in(\Omega \times \bfR \times \bfR^{n})$, then
$u\leq v$ in $\Omega$.
\end{lemma}
\begin{proof}
Suppose the contrary that for some $y\in \Omega$,
$$u(y)-v(y)=\max_{\Omega}(u-v)>0.$$
Then $S_{k}(D^2v(y))\geq S_{k}(D^2u(y)),$ as the Hessian
$D^2(v-u)$ is positive semi-definite at $y$ and $u,v$ are
$k$-convex. On the other hand ,we use (2.1)-(2.4) and the facts
$u(y)>v(y)$ and $Du(y)=Dv(y)$ to obtain
$$ S_{k}(D^2v(y))\leq \phi(y,v(y),Dv(y))<\psi(y,u(y),Du(y))\leq S_{k}(D^2u(y)).$$
This is a contradiction.
\end{proof}
\begin{remark}
The assumption we have used is $\psi_{z}(y,u(y),Du(y))>0 $ (or
$\phi_{z}(y,u(y),Du(y))>0 $) at the point $y$ instead of
$\psi_{z}(x,z,\vp)>0 $ (or $\phi_{z}(x,z,\vp)>0 $ ) for all $(x,
z,\vp)\in(\Omega \times \bfR \times \bfR^{n})$.
\end{remark}
\begin{theorem}
Assume $u=+\infty, v=+\infty$ on $\partial\Omega$, $u$ is
$k$-convex and $v$ is strictly $k$-convex in $\Omega ,$ satisfying
(2.1) and (2.2).
 Suppose the domain $\Omega$ is bounded and
star-shaped with respect to a point $x_{0}$ and $\psi$ satisfies
\begin{equation}
x\cdot D_{x}\psi(x,z,\vp)\leq 0 \; \mbox{and $\vp\cdot
D_{\vp}\psi(x,z,\vp)\geq 0$}, \;\;\forall
\;(x,z,\vp)\in\Omega\times\bfR\times\bfR^{n}.
\end{equation}
If, in addition, either there is a $q>k$ such that
\begin{equation}
\psi(x,\mu z^{+},\vp)\geq\mu^q \psi(x,z,\vp), \;\; \forall \;
\mu\geq 1, \forall \;(x, z,\vp)\in(\Omega \times \bfR \times
\bfR^{n})
\end{equation}
or there is a $\varepsilon>0$ such that
\begin{equation}
\psi_{z}(x,z,\vp)\geq\varepsilon \psi(x,z,\vp), \;\; \forall \;(x,
z,\vp)\in(\Omega \times \bfR \times \bfR^{n}),
\end{equation}
then $u\leq v$ in $\Omega$. Particularly, problem (1.1)-(1.2) has
at most one strictly $k$-convex solution in $C^{2}(\Omega)$.
\end{theorem}
\begin{proof}
Without loss of generality we may assume $x_{0}=0$. For
$\lambda\in(0,1)$, let
$$u_{\lambda}(x)=\lambda^{\alpha}u(\lambda x)-a,\;\;x\in\Omega_{\lambda}$$
where $\Omega_{\lambda}=\{x\in\bfR^{n}:\lambda x\in \Omega\}$ and
\[ \begin{cases}
  a = 0,  \;\alpha = \frac{2k}{q-k}, & \mbox{if (2.6) holds}, \\
  a=-\frac{2k}{\varepsilon}\ln\lambda, \;\alpha = 0, &  \mbox{if (2.7) holds}.
     \end{cases} \]
  Using (2.1)
and (2.5), we have
\begin{equation}
 \begin{aligned}
 \;\;S_{k}(D^2 u_{\lambda}(x))
& = \lambda^{k(2 + \alpha)}S_{k}(D^2u(\lambda x)) \\
& \geq \lambda^{k(2 + \alpha)}
       \psi\big(\lambda x, u (\lambda x), Du (\lambda x)\big) \\
& = \lambda^{k(2 + \alpha)} \psi \big(\lambda x,
       \lambda^{- \alpha} (u_{\lambda} (x) + a),
       \lambda^{- (1+ \alpha)} Du_{\lambda} (x)\big) \\
& = \lambda^{k(2 + \alpha)} \int^{1}_{0}\partial_{t}\psi
\big(x+t(\lambda
x-x),\lambda^{- \alpha} (u_{\lambda} (x) + a),\\
&Du_{\lambda} (x)+t(\lambda^{- (1+ \alpha)} Du_{\lambda}
(x)-Du_{\lambda} (x))\big)dt\\
&+ \lambda^{k(2 + \alpha)}\psi \big(x,
       \lambda^{- \alpha} (u_{\lambda} (x) + a), Du_{\lambda} (x)\big)\\
& \geq\lambda^{k(2 + \alpha)}\psi \big(x,
       \lambda^{- \alpha} (u_{\lambda} (x) + a), Du_{\lambda}
       (x)\big).
 \end{aligned}
\end{equation}
Note that (2.6) implies $\psi\leq0$ where $z\leq0$ and $\psi=0$ at
$z=0.$ We conclude that
\begin{equation}
S_{k}(D^2 u_{\lambda}(x)) \geq\psi \big(x, u_{\lambda} (x),
Du_{\lambda}(x)\big), \forall x\in \Omega .
\end{equation} In fact, if  (2.6) holds,  then $a=0$ and $\alpha=\frac{2k}{q-k}.$ Hence,
 we use (2.8) and (2.6) to obtain
$$ S_{k}(D^2 u_{\lambda}(x))\geq\lambda^{k(2 + \alpha)-\alpha q}\psi \big(x, u_{\lambda} (x),
Du_{\lambda}(x)\big)=\psi \big(x, u_{\lambda} (x),
Du_{\lambda}(x)\big).$$  Note  that (2.7) implies
$$ \psi(x,z_{1},\vp)\geq e^{\varepsilon(z_{1}-z_{2})}\psi(x,z_{2},\vp), \ \ \forall z_1, z_2\in R.$$
So, if   (2.7) holds,  then $a=-\frac{2k}{\varepsilon} $ and
$\alpha =0 .$ Hence, we use (2.8) to obtain
$$S_{k}(D^2 u_{\lambda}(x))\geq\lambda^{2k}e^{\varepsilon a}\psi \big(x, u_{\lambda} (x),
Du_{\lambda}(x)\big)=\psi \big(x, u_{\lambda} (x),
Du_{\lambda}(x)\big). $$ Therefore, (2.9) holds in any case. We
claim that
\begin{equation}
v\geq u_{\lambda}\;\;\mbox{on $\Omega$}\;\;\mbox{for all
$\lambda\in (0,1)$}.
\end{equation}
Suppose the contrary that $u_{\lambda}(x_{0})>v(x_{0})$ for some
$x_{0}\in \Omega$. Since
$\overline{\Omega}\subset\Omega_{\lambda}$ and
$v-u_{\lambda}=+\infty$ on $\partial\Omega$, we have a $y\in
\Omega$, such that
$$u_{\lambda}(y)-v(y)=\max_{\Omega}(u_{\lambda}-v)>0  .$$
Hence, by(2.5),(2.2),(2.1) and the strict $k$-convexity of $v$,
\[ \begin{aligned}
  \psi(y,u_{\lambda}(y),Du_{\lambda}(y))
&\geq\psi(y,v(y),Du_{\lambda}(y))\\
& =\psi(y,v(y),Dv(y))\\
& \geq\phi(y,v(y),Dv(y))\\
& \geq S_{k}(D^{2}v)\\
&>0
 \end{aligned} \]
which, together with (2.6) or (2.7), implies
$\psi_{z}(y,u_{\lambda}(y),Du_{\lambda}(y))>0$. (Note that if
(2.6) holds, then (2.2) implies $\psi\leq0$ where $z\leq0$. So
$v(y)>0$  and $u_{\lambda}(y)>0$ then (2.6) implies
$\psi_{z}(y,u_{\lambda}(y),Du_{\lambda}(y))>0$). Consequently, we
obtain a contradiction as in the proof of Lemma 2.1( See Remark
2.2). This proves (2.10). Letting $\lambda\rightarrow 1^{-}$ in
(2.10), we obtain the desired result.
\end{proof}

\bigskip

\section{Barriers and non-existence}
\setcounter{equation}{0}

In this section, we construct some barriers that will be used in
the proof of our main results. In particular, we will prove
theorems 1.1 and 1.2.

Let $ u(x)=u(|x|) $ be a radially symmetric function. A
straightforward calculation gives
\begin{equation}
S_{k}(D^{2}(u))=A_{k}r^{1-n}[r^{n-k}(u^{\prime})^{k}]^{\prime}
,\;\; r=|x|,
\end{equation}
where $ A_{k}=\frac{(n-1)!}{k!(n-k)!} $ (See[15; p.285]). Hence
equation (1.1) is written as
\begin{equation}
A_{k}[(r^{\frac{n}{k}-1}u^{\prime})^{k}]^{\prime}=r^{n-1}
\psi(x,u,\frac{xu^{\prime}}{r}).
\end{equation}
\begin{lemma}
Let $\eta\in C^1 (\bfR)$ satisfy $\eta(z)>0,\;
\eta^{\prime}(z)\geq0$ for all $z\in \bfR.$ Then, for any $a>0,$
there exists a strictly convex radially symmetric function $v\in
C^{2}(B_{a}(0)) $ satisfying
\begin{equation}
 \begin{cases} S_{k}(D^{2}v) \geq e^{v}\eta(v)(1+{|Dv|}^{k}) \; \mbox{in
 $B_{a}(0)$},\\
 v(0) \leq 0 ,\; v = +\infty  \; \mbox{ on $\partial
 B_{a}(0)$}.
 \end{cases}
\end{equation}
\end{lemma}
\begin{proof}
Consider the initial value problem
\begin{equation}
 \begin{cases}
 \varphi^{\prime}=[exp(A^{-1}_{k}r^{k}e^{\varphi}\eta(\varphi))-1]^{\frac{1}{k}},
 \;\; r>0\\ \varphi(0)=0.
 \end{cases}
\end{equation}
Let $[0,T)$ be the maximal interval on which the solution to (3.4)
exists. We conclude that $T$ is finite. In fact, it follows from
(3.4) and $\eta^{\prime}\geq 0$ that
$$\varphi^{\prime}(r)\geq
r[A^{-1}_{k}e^{\varphi}\eta(\varphi)]^{\frac{1}{k}}\geq
r[A^{-1}_{k}e^{\varphi(r)}\eta(0)]^{\frac{1}{k}},\;0<r<T.
$$
Since $\varphi(0)=0,$ we have
$$
k\geq
k(1-e^{\frac{-\varphi(\rho)}{k}})=\int^{\rho}_{0}\varphi^{\prime}(r)e^{\frac{-\varphi(r)}{k}}dr\geq
\big(\frac{\eta(0)}{A_{k}}\big)^{\frac{1}{k}}\int^{\rho}_{0}rdr
=\frac{1}{2}\big(\frac{\eta(0)}{A_{k}}\big)^{\frac{1}{k}}\rho^{2}
$$
for any $\rho\in(0,T)$. This proves $T<\infty.$ Furthermore, we
see that $\varphi \in C^{2}[0,T)$ and $\varphi(T)=+\infty$ as
$\varphi^{\prime}(\rho)>0$ for $\rho> 0.$ It follows from (3.4)
that
$$\ln(1+(\varphi^{\prime})^{k})=A_{k}^{-1}r^{k}e^{\varphi}\eta(\varphi),$$
whose differentiation yields
\begin{equation}
\frac{k(\varphi^{\prime})^{k-1}\varphi^{\prime\prime}}{1+(\varphi^{\prime})^{k}}
\geq kA_{k}^{-1}r^{k-1}e^{\varphi}\eta(\varphi).
\end{equation}
This, particularly, implies $\varphi^{\prime\prime}> 0$ in
$(0,T).$

For given $a> 0$, define $v$ by
$$v(x)=\varphi\big(\frac{T|x|}{a} \big)-2k\big(-\ln\frac{T}{a} \big)^{+},\; x\in B_{a}(0).$$
Then $v\in C^{2}(B_{a}(0)),\; v(0) \leq 0 ,\; v = +\infty  \;
\mbox{on $\partial B_{a}(0)$}$, and it is strictly convex, since
$\varphi\in C^{2}[0,T)$ and $\varphi^{\prime\prime}> 0.$ Using
(3.1) and (3.5), we compute for $x\in B_{a}(0)$ that
\[ \begin{aligned}
S_{k}(D^2 v (x))
&=\big(\frac{T}{a}\big)^{2k}A_{k}\big(\frac{T}{a}|x|\big)^{1-n}\big[k\big(\frac{T}{a}|x|\big)^{n-k}
\big(\varphi^{\prime}\big(\frac{T}{a}|x|\big)\big)^{k-1}\varphi^{\prime\prime}\big(\frac{T}{a}|x|\big)\\
& +(n-k)\big(\frac{T}{a}|x|\big)^{n-k-1}\big(\varphi^{\prime}\big(\frac{T}{a}|x|\big)\big)^{k}\big]\\
& \geq
\big(\frac{T}{a}\big)^{2k}A_{k}\big(\frac{T}{a}|x|\big)^{1-k}k\big(\varphi^{\prime}\big(\frac{T}{a}|x|\big)\big)^{k-1}
\varphi^{\prime\prime}\big(\frac{T}{a}|x|\big)  \\
\medskip
& \geq
\big(\frac{T}{a}\big)^{2k}e^{\varphi\big(\frac{T}{a}|x|\big)
}\eta \big( \varphi\big(\frac{T}{a}|x|\big)\big)\big[1+ \big(\varphi^{\prime}\big(\frac{T}{a}|x|\big)\big)^{k}\big]\\
\medskip
& \geq \big(\frac{T}{a}\big)^{2k}e^{v(x)+2k\big(-\ln\frac{T}{a}
\big)^{+}}\eta\big(v(x)+2k\big(-\ln\frac{T}{a} \big)^{+}\big)\big[1+\big(\frac{T}{a}\big)^{-k}|Dv(x)|^{k}\big]\\
\medskip
& \geq e^{v (x)} \eta \big(v (x)\big) \big(1 + |Dv (x)|^k\big).
\end{aligned}  \]
\end{proof}
\begin{notation}
Given $a$ and $\eta$, we will use $v^{a,\eta}(x)=v^{a,\eta}(|x|)$
to denote the function $v\in C^{2}(B_{a}(0))$ obtained as in Lemma
3.1.
\end{notation}
\begin{lemma}
Let $\Omega$ be a domain contained in a ball $B_{a}(x_{0})$ and
$u\in C^{2}(\Omega)$ a $k$-convex solution of (1.1)-(1.2). If
there exists a function $\eta\in C^{1}(\bfR),\; \eta>0$ and
$\eta^{\prime}\geq 0$ in $\bfR$, such that
$$\psi(x,z,\vp)\leq e^{z}\eta(z)(1+|\vp|^{k}),\;\;\forall\;(x,z,\vp)\in \overline{\Omega}\times\bfR\times\bfR^{n},$$
then $u(x)\geq v^{a,\eta}(x-x_{0})$ for all $x\in \Omega.$
\end{lemma}
\begin{proof}
Without loss of generality, we assume $x_{0}=0.$ For any $r>a$,
since $u-v^{r,\eta}=+\infty$ on $\partial\Omega$, by Lemma 2.1 and
Lemma 3.1, we have $u\geq v^{r,\eta}$ in $\Omega.$ Letting
$r\rightarrow a^{+}$, we obtain $u\geq v^{a,\eta}.$
\end{proof}
Next for $q>k$, consider the function
\begin{equation}
w(x)=(1-|x|^{2})^{\frac{k+1}{k-q}}=(1-r^{2})^{\frac{k+1}{k-q}}=w(r),\;\;r=|x|.
\end{equation}
Observing that $w^{\prime}\geq 0$ and $w^{\prime\prime}>0,$ by
(3.1) we have for any $r\in [0,1)$ that
\[ \begin{aligned}
S_{k}(D^{2}w(x))
&=A_{k}[(n-k)r^{-k}(w^{\prime})^{k}+kr^{1-k}(w^{\prime})^{k-1}w^{\prime\prime}]\\
&\leq
C_{1}(n,k,q)[(1-r^{2})^{(\frac{q+1}{k-q})k}+(1-r^{2})^{(\frac{q+1}{k-q})(k-1)}(1-r^{2})^{\frac{q+1}{k-q}-1}]\\
&=C_{2}(n,k,q)(1-r^{2})^{(\frac{q+1}{k-q})k-1}[(1-r^{2})+1]\\
&\leq 2C_{2}(n,k,q)(1-r^{2})^{(\frac{k+1}{k-q})q}\\
&=2C_{2}(n,k,q)w^{q}(x).
\end{aligned}  \]
Hence, we have a positive constant $B=B(n,k,q)$ such that
\begin{equation}
S_{k}(D^{2}w)\leq Bw^{q} \;\;\mbox{in $B_{1}(0)$}.
\end{equation}
By rescaling, we define
\begin{equation}
w^{a,M}(x)=\lambda w(\frac{x}{a}),\;\; x\in B_{a}(0),\;\;
\lambda=\big(\frac{B}{a^{2k}M}\big)^{\frac{1}{q-k}}.
\end{equation}
\begin{lemma}
For any $a,M> 0$ and $q>k,$ $ w^{a,M}\in
C^{\infty}(B_{a}(0)),\;w^{a,M}=+\infty$ on $\partial B_{a}(0)$ and
$$ S_{k}\big(D^{2}w^{a,M}\big)\leq M(w^{a,M})^{q} \;\;\mbox{in $ B_{a}(0)$}. $$
\end{lemma}
\begin{proof}
By a direct calculation, (3.7) and (3.8), we have
$$S_{k}\big(D^{2}w^{a,M}(x)\big)=\frac{\lambda^{k}}{a^{2k}}S_{k}\big(D^{2}w(\frac{x}{a})\big)
\leq\frac{\lambda^{k}B}{a^{2k}}w^{q}\big(\frac{x}{a}\big)=M\big(w^{a,M}(x)\big)^{q}\;\;\mbox{in
$B_{a}(x)$}.$$
\end{proof}
\begin{lemma}
Let $u\in C^{2}(\Omega)$ be a $k$-convex solution of (1.1). If
$\Omega$ contains a ball $B_{a}(x_{0})$ and $\psi$ satisfies (1.8)
for some $q>k$ and $M>0,$ then $u(x)\leq w^{a,M}(x-x_{0})$ in
$B_{a}(x_{0})$.
\end{lemma}
\begin{proof}
It is immediate from Lemmas 3.4 and 2.1.
\end{proof}

 Note that for any domain $\Omega$ and any $x\in \Omega$, the ball
 $B_{d(x)}(x)\subset\Omega$, where $d(x)$ is the distance function to $\partial
 \Omega$. Then Lemma 3.5 implies
\begin{corollary}
 Let $u \in C^2 (\Omega)$ be a $k$-convex solution of (1.1). If $\psi$ satisfies (1.8) for some $q > k$
and $M > 0$, then
\[ u (x) \leq \bar{h} (d (x)), \;\; \forall \; x \in \Omega \]
where $\bar{h} \in C^{\infty} (\bfR^+)$ is given by
\begin{equation}
 \bar{h} (r) = w^{r, M} (0), \;\; r > 0.
\end{equation}
\end{corollary}
\begin{lemma}
Assume $\gamma, q \geq 0$, $ \gamma+q \leq k$ and $M > 0$. Then
there exists a  strictly convex radially symmetric positive
function $\tilde{u} \in C^{\infty} (\bfR^n)$ satisfying
\begin{equation}
 S_{k}(D^2 \tilde{u} (x)) \geq
    M \big(1 + (\tilde{u} (x))^q\big) \big(1 + |D\tilde{u} (x)|^\gamma\big),
\;\; \forall \; x \in \bfR^n.
\end{equation}
\end{lemma}
\begin{proof}
We want only to combine the following three conclusions. First,
for any $p^{\prime},q^{\prime}\geq 0, p^{\prime}+q^{\prime}\leq
n,$ and $M^{\prime}>0,$ by Lemma 3.7 in [7], one has a strictly
convex radially symmetric positive function $\tilde{u} \in
C^{\infty} (\bfR^n)$ satisfying
\begin{equation}
 S_{n}(D^2 \tilde{u} (x)) \geq
    M^{\prime} \big(1 + (\tilde{u} (x))^{p^{\prime}}\big) \big(1 + |D\tilde{u} (x)|^{q^{\prime}}\big),
\;\; \forall \; x \in \bfR^n.
\end{equation}
Second, for each $k\in \{1,2,...,n-1 \}$ and $\lambda\in
\Gamma_{k+1},$ it follow from [10] or [18] that
$$\frac{(k+1)!(n-k-1)!}{n!}S_{k+1}(\lambda)\leq \big[\frac{k!(n-k)!}{n!}S_{k}(\lambda)\big]^{\frac{k+1}{k}}.$$
Hence, there is a positive constant $C_{1}=C(n,k)$ such that
\begin{equation}
C_{1}S_{n}^{\frac{k}{n}}(\lambda)\leq
S_{k}(\lambda),\;\;\forall\;\lambda\in \Gamma_{n}.
\end{equation}
Finally, we can choose a positive constant $C_{2}=C(n,k)$ such
that
\begin{equation}
(1+t)^{\frac{k}{n}}\geq C_{2}(1+t ^{\frac{k}{n}}),\;\;\forall
\;t\geq 0.
\end{equation}
Combing this with (3.11) and (3.12), we obtain (3.10).
\end{proof}
\begin{proof}[Proof of Theorem 1.1] Let $u \in C^2 (\Omega)$ be a
$k$-convex solution of (1.1)-(1.2). We will induce a
contradiction.

 Let $\tilde{u} $ be the same function as in Lemma 3.7, where $\gamma,q$ and $M$ are as in (1.5).
 Observing that $u - C\tilde{u} = +\infty$ on $\partial \Omega$ for any $C >
 0$, we can choose $C>1$ and a $y\in \Omega$ such that
\[ u (y)- C \tilde{u} (y) = \min_{\Omega} (u - C \tilde{u}) < 0.
\]
Hence $Du(y)=CD\tilde{u}(y)$ and
$\big(D^{2}u(y)-CD^{2}\tilde{u}(y)\big)$ is a positive
semi-definite matrix. However, it follows from (1.1) and (1.5)
that
\[ \begin{aligned}
S_{k}(D^2 u (y))
   & \leq M (1 + \big(u^+ (y))^q\big) \big(1 + |Du (y)|^\gamma\big) \\
   & < M \big(1 + (C \tilde{u} (y))^q\big)
                  \big(1 + C|D \tilde{u} (y)|^\gamma\big) \\
   &  <   C^k M \big(1 + (\tilde{u} (y))^q\big)
                  \big(1 + |D \tilde{u} (y)|^\gamma\big) \\
   & \leq C^k S_{k}( D^2 \tilde{u} (y))\\
   & = S_{k} (CD^2 \tilde{u} (y)),
   \end{aligned} \]
 a contradiction.
\end{proof}

In order to prove Theorem 1.2, we need the following
\begin{lemma}
For any $\alpha>1$ and $a>0,$ there exists a strictly convex
radially symmetric function $\bar{u} \in C^2 (B_a (0))$ satisfying
$$S_{k}(D^2 \bar{u}) \leq \frac{(k+1)(n-1)!}{a^k
  (\alpha-1)k!(n-k-1)!} \big(1 + |D\bar{u}|^k\big)^{\alpha}
        \;\; \mbox{in $B_a (0)$}$$
and
$$ \frac{\partial \bar u}{\partial \nu} = + \infty  \;\;
\mbox{on $\partial B_a (0)$},$$ where $\nu$ is the unit normal to
$\partial B_a (0).$
\end{lemma}
\begin{proof}
Let $ \beta=\frac{1}{\alpha-1}$ and
\[\varphi(r)=
\begin{cases}
\int_{0}^{r}[\frac{(1-t^{k+1})^{-\beta}-1}{t}]^{\frac{1}{k}}dt,\;\;
& r\in (0,1), \\
\;\;\;\;\;\;\;\;\;\;\;\;\;\;0, & r=0.
\end{cases}\]
It is easy to verify that
\begin{equation}
\varphi\in C^{2}[0,1),\;\; \varphi(0)=\varphi^{\prime}(0)=0; \;\;
1+r(\varphi^{\prime}(r))^{k}=(1-r^{k+1})^{-\beta},\;\; \forall
\;0\leq r<1.
\end{equation}
and
\begin{equation}
\begin{aligned}
(\varphi^{\prime}(r))^{k}+kr(\varphi^{\prime}(r))^{k-1}\varphi^{\prime\prime}(r)
&=\frac{k+1}{\alpha-1}r^{k}\big(1+r(\varphi^{\prime}(r))^{k}\big)^{\alpha}\\
&=(k+1)\beta r^{k}(1-r^{k+1})^{-\beta-1},\;\;\forall \;r\in(0,1).
\end{aligned}
\end{equation}
We claim that
\begin{equation}
\varphi^{\prime}\geq 0 \;\; \mbox{and}\;\;
\varphi^{\prime\prime}>0 \;\;\mbox{in}\;[0,1).
\end{equation}
In fact, a direct differentiation yields
$$\varphi^{\prime}(r)=\big[\frac{(1-r^{k+1})^{-\beta}-1}{r}\big]^{\frac{1}{k}},\;\; \forall \;r\in(0,1).$$
Then
\begin{equation}
\lim_{r\rightarrow0^{+}}\varphi^{\prime}(r)=0,\;\;\lim_{r\rightarrow1^{-}}\varphi^{\prime}(r)=+\infty,\;\;
\varphi^{\prime}\geq 0 \;\;\mbox{in}\;[0,1)
\end{equation}
and
\[ \begin{aligned}
\varphi^{\prime\prime}(0)
&=\lim_{r\rightarrow0^{+}}\frac{\big[(1-r^{k+1})^{-\beta}-1\big]^{\frac{1}{k}}}{r^{\frac{1}{k}+1}}\;\;
(\mbox{letting $t=r^{k+1}$})\\
&=\lim_{t\rightarrow0^{+}}\big[\frac{(1-t)^{-\beta}-1}{t}\big]^{\frac{1}{k}}\\
&=\beta^{\frac{1}{k}}>0.
\end{aligned} \]
For $r\in(0,1),$ it follows from (3.15) that
$\varphi^{\prime\prime}(r)$ has the same sign as $g_{1}(r^{k+1})$
where
$$g_{1}(t)=\beta(k+1)t(1-t)^{-\beta-1}-(1-t)^{-\beta}+1,\;\;\forall\;t\in(0,1).
$$
Let $1-t=s.$ We see that $g_{1}(t)=g_{2}(s)$ where
$$g_{2}(s)=\beta(k+1)(1-s)s^{-\beta-1}-s^{-\beta}+1,\;\;\forall\;s\in(0,1).$$
Furthermore, $g_{2}(s)$ has the same sign as $g_{3}(s)$ where
$$g_{3}(s)=\beta(k+1)-\beta(k+1)s-s+s^{\beta+1},\;\;\forall\;s\in(0,1).$$
Since $g_{3}^{\prime}(s)<0$ for $s\in(0,1)$ and $g_{3}(1)=0,$ we
see that $g_{3}(s)>0$ for all $s\in(0,1).$ Therefore, we obtain
(3.16).

Now, let $\bar{u}(x)=a\varphi(a^{-1}|x|),\;\;x\in B_{a}(0).$ Then
$ \bar{u}\in C^{2}(B_{a}(0))$ and is strictly convex by (3.16),
satisfying
 $\frac{\partial \bar u}{\partial \nu} =+\infty $ on $\partial
B_a (0)$ by (3.17). Moreover by (3.1),(3.15) and (3.16), we have
\[ \begin{aligned}
S_{k}(D^2\bar{u}(x))
&=\frac{A_{k}}{a^{k}}\big(\frac{|x|}{a}\big)^{-k}
\big[(n-k)(\varphi^{\prime}(\frac{|x|}{a}))^{k}+k\frac{|x|}{a}
(\varphi^{\prime}(\frac{|x|}{a}))^{k-1}\varphi^{\prime\prime}(\frac{|x|}{a})\big]\\
&\leq \frac{(n-k)(k+1)A_{k}}{a^{k}(\alpha-1)}
\big[1+\frac{|x|}{a}(\varphi^{\prime}(\frac{|x|}{a}))^{k}\big]^{\alpha}\\
&\leq\frac{(n-k)(k+1)A_{k}}{a^{k}(\alpha-1)}
\big[1+|D\bar{u}(x)|^{k}\big]^{\alpha} \\
&=\frac{(k+1)(n-1)!}{a^{k}(\alpha-1)k!(n-k-1)!}
\big[1+|D\bar{u}(x)|^{k}\big]^{\alpha}\;\;\mbox{in}\;B_a (0).
\end{aligned}  \]
\end{proof}
\begin{proof}[Proof of Theorem 1.2]
We may assume $\Omega\supset\overline{B_{a}(0)}.$ Suppose the
contrary that there was a $k$-convex solution $u\in C^{2}(\Omega)$
to (1.1)-(1.2). We will derive a contradiction. Let $a, \alpha$
and $M$ be the same as in Theorem 1.2. Choose a function $\bar{u}$
as in Lemma 3.8. Then we have
\begin{equation}
S_{k}(D^{2}\bar{u}(x))<M\big(1+|D\bar{u}(x)|^{k}\big)^{\alpha}\;\;\mbox{in}\;B_{a}(0),
\end{equation}
and
$$\frac{\partial}{\partial \nu}(\bar{u}-u)
= + \infty  \;\; \mbox{on $\partial B_a (0)$}.$$ It follows from
the last equation that
$$\bar{u}(y)-u(y)=\min_{B_a (0)}(\bar u-u) $$
for some $ y\in B_a (0).$ Using (1.1),(1.6) and (3.18), and
repeating the same arguments as in the proof of Theorem 1.1, we
obtain a contradiction immediately.
\end{proof}
\bigskip

\section{Proof of theorem 1.3}
\setcounter{equation}{0} We divide the proof into two steps.
\begin{proof}   {\it Step 1.}
Assume $\Omega$ is a bounded strictly convex smooth domain. We
will find a solution of (1.1)-(1.2) as required as in Theorem 1.3
by the limit of solutions, $u_{m}$, of the following Dirichlet
problem
\begin{equation}
\begin{cases}
S_{k}(D^{2}u)=\psi(x,u,Du) & \mbox{in $\Omega$}\\
 \;\;\;\;\;\;\;\;\;\;\;u=m & \mbox{on $ \partial\Omega $}
\end{cases}
\end{equation}
where $m=1,2,3,...$ . By assumption (1.10), we may find a positive
nondecreasing function $\eta\in C^{\infty}(\bfR^{n})$ such that
\begin{equation}
\max_{y\leq z}\phi(y)\leq e^{\varepsilon z}\eta(z),\;\;\forall
\;z\in \bfR.
\end{equation}
Without loss of generality, we assume $\varepsilon=1$ as this may
be achieved by rescaling. Since $\Omega$ is bounded, we may choose
$a>0$ such that $\Omega\subset B_{a}(0)$ and $v^{a,\eta}\leq 1$ on
$\partial\Omega$ (See (3.3) and Notation 3.2). Using (4.2), (1.9),
(1.10) and applying (1.7) and Lemma 2.1, we obtain that for any
$k$-convex solution $u_{m}\in C^{2}(\overline{\Omega})$ to (4.1),
\begin{equation}
v^{a,\eta}\leq u_{m}(x)\leq m, \;\;\forall
\;x\in\Omega,\;\;\forall \; m\geq 1.
\end{equation}
Let
$$C^{0}_{m}=\max\{m,\sup_{\Omega}|v^{a,\eta}|\},\;\;\forall \; m \geq 1.$$
In order to show the existence of (4.1), we want to use results of
Lions [12] as well as of Guan [6]. First, by a result of [12],
there exists, for each $m$ and any constant $C_{m}>0$, a strictly
convex function $\underline{u}_{m}\in C^{2}(\overline{\Omega})$
satisfying
\begin{equation}
\begin{cases}
\det(D^{2}\underline{u}_{m})\geq C_{m}(1+|D\underline{u}_{m}|^{n}) & \mbox{in $\Omega$}\\
 \;\;\;\;\;\;\;\;\;\;\;\;\;\underline{u}_{m}=m & \mbox{on $ \partial\Omega
 $}.
\end{cases}
\end{equation}
Using this, (3.12) and (3.13), and choosing a suitable $C_{m}$, we
see that
\begin{equation}
\begin{cases}
S_{k}(D^{2}\underline{u}_{m})\geq \phi(C^0_{m})(1+|D\underline{u}_{m}|^{k}) & \mbox{in $\Omega$}\\
 \;\;\;\;\;\;\;\;\;\;\;\;\;\underline{u}_{m}=m & \mbox{on $ \partial\Omega
 $},
\end{cases}
\end{equation}
which means that $\underline{u}_{m}$ is a subsolution of (4.1) for
each $m$, since $\phi $ is nondecreasing and $m\geq
\underline{u}_{m}$ in $\Omega.$  This fact, together with (1.7),
(1.9), (1.11),  Theorem 1.2 of Guan [6], implies that problem
(4.1) has a unique $k$-convex solution $u_{m}\in
C^{\infty}(\overline{\Omega})$ for each $m.$ Moreover, we have
\begin{equation}
u_{m}(x)\leq u_{m+1}(x),\;\;\forall \; x\in\Omega,\;\;\forall
\;m\geq 1
\end{equation}
 by Lemma 2.1. We claim that there exists $a>0$ depending only on
$\Omega$ and an decreasing sequence $a_{m}\rightarrow a
(m\rightarrow \infty)$ such that
\begin{equation}
v^{a_{m},\eta}(a-d(x))\leq u_{m}(x)\leq \bar{h}(d(x)),\;\;
\forall\;x\in\Omega,\;\forall m\geq1.
\end{equation}
In fact, the second inequality in (4.7) follows directly from
Corollary (3.6). To show the first one, we use the strict
convexity of $\Omega$ to find the smallest positive number $a$,
such that for any $\bar{x}\in \partial\Omega,$ there is a ball
$B_{a}(x_{0})\supset\Omega$ satisfying
$\overline{\Omega}\cap\partial B_{a}(x_{0})=\{\bar{x}\}.$ Choose
$a_{1}>a_{2}>...>a_{m}>a_{m+1}>...$ , $a_{m}\rightarrow
a(m\rightarrow \infty),$ such that $v^{a_{m},\eta}(a)=m$ for each
$m\geq1.$ For any $y\in \Omega,$ choose $\bar{y}\in
\partial\Omega$ and a ball $B_{a}(x_{0})$ such that
$d(y)=|y-\bar{y}|$, $\Omega\subset B_{a}(x_{0}) $,
$\overline{\Omega}\cap\partial B_{a}(x_{0})=\{\bar{y}\}.$
Observing that $$v^{a_{m},\eta}(x-x_{0})\leq
v^{a_{m},\eta}(a)=m\leq u_{m}(x),\;\;\forall \;x\in
\partial\Omega,$$
we use  Lemma 2.1 to obtain
$$v^{a_{m},\eta}(x-x_{0})\leq u_{m}(x),\;\;\forall \;x\in\Omega.$$
In particular,
$$v^{a_{m},\eta}(a-d(y))=v^{a_{m},\eta}(y-x_{0})\leq u_{m}(y).$$
This proves the first inequality in (4.7), since $y\in\Omega$ is
arbitrary.

Now by (4.6) and (4.7), we see that for each $x\in\Omega$, the
limit $$u(x)= \lim_{m\rightarrow\infty} u_{m}(x)$$ exists and it
satisfies
\begin{equation}
v^{a,\eta}(a-d(x))\leq u(x)\leq \bar{h}(d(x)),\;\;
\forall\;x\in\Omega.
\end{equation}
Moreover, by Theorem 3.1 of [17], the convergence is uniform in
every compact set $K\subset\Omega$ and $u\in
C_{loc}^{0,1}(\Omega).$ By the stability theorem of viscosity
solutions under the uniform convergence, we see that $u$ is a
viscosity $k$-convex solution of (1.1)-(1.2), which satisfies
(1.12) by (4.8). \vskip 0.5cm
 {\it Step 2.}
Suppose now that $\Omega$ is a bounded strictly convex domain. We
will complete the proof of Theorem 1.3 . In this case , we choose
a sequence of strictly convex smooth domains
\begin{equation}
\Omega_{1}\subset\Omega_{2}\subset...\subset
\Omega_{m}\subset\Omega_{m+1}\subset...\subset\Omega,
\end{equation}
such that
$$\Omega=\bigcup_{m=1}^{\infty}\Omega_{m}.$$
For each $m\geq 1,$ by the result of Step 1, we choose a
$k$-convex viscosity solution $u_{m}\in C_{loc}^{0,1}(\Omega _m)$
to the problem
\[\begin{cases}
S_{k}(D^{2}u)=\psi(x,u,Du) & \mbox{in $\Omega_{m}$},\\
 \;\;\;\;\;\;\;\;\;\;\;u=+\infty & \mbox{on $ \partial\Omega_{m} $}.
\end{cases}
\]
By (4.8), (4.9) and Lemma 2.1, we may assume
\begin{equation}
\begin{aligned}
v^{a_{m+1},\eta}(a_{m+1}-d_{m+1}(x))&\leq u_{m+1}(x)\\
 &\leq
u_{m}(x) \leq \bar{h}(d_{m}(x)),\;\; \forall\;x\in\Omega_{m},\;\;
\forall \;m\geq1,
\end{aligned}
\end{equation}
where $d_{m}(x)=dist(x,\partial\Omega_{m})$ and $a_{m}$ is the
smallest positive number, such that for any
$\bar{x}\in\partial\Omega_{m},$ there is a ball
$B_{a_{m}}(x_{0})\supset\Omega_m$ satisfying
$\overline{\Omega}_m\cap
\partial B_{a_{m}}(x_{0})=\{\bar{x}\}.$
Note that (4.9)  implies $\{a_{m}\}$ is a   nondecreasing
sequence. Furthermore, $\{a_{m}\}$ is bounded since $\Omega$ is a
bounded strictly convex. Letting
 $a=\lim_{m\rightarrow\infty}a_{m}$ and using (4.10), we see
 that the limit function $$u(x)=\lim_{m\rightarrow\infty}u_{m}(x)$$
 exists for each $x\in \Omega$   and it satisfies (4.8) again.
 Repeating the arguments from (4.8) to the end in Step 1, we have
 completed the proof of Theorem 1.3 .
\end{proof}

\vskip 0.5cm

{\bf Acknowledgement.}  The author would like to thank the referee
and Professor Bo Guan for useful comments and suggestions.  Dr.
Xiuqing Chen read throughout this paper and typed it. The author
thanks for his work and suggestions.


\begin{thebibliography}{99}

\bibitem{1}
J. Bao, J. Chen, B. Guan and M. Ji, {\em Liouville property and
regularity of a Hessian quotient equation}, {Amer. J. Math.},{\bf
125}(2003), 310-316.

\bibitem{ 2}
L. A. Caffarelli, L. Nirenberg, J. Spruck, {\em The Dirichlet
problem for nonlinear second order elliptic
equations,III:Functions of the eigenvalues of the Hessian}, {Acta
Math.}, {\bf 155} (1985), 261-301.

\bibitem{3}
S. Y. Cheng and S. T. Yau, {\em On the existence of a complete
K\"ahler metric on non-compact complex manifolds and the
regularity of Fefferman's equation}, { Comm. Pure Applied Math.},
{\bf 33} (1980), 507-544.

\bibitem{4}
S. Y. Cheng and S. T. Yau, {\em The real Monge-Amp\`ere equation
and affine flat structure}, Proc. 1980 Beijing Symp. on Diff.
Geom. and Diff. Equations Vol I, 339-370 (1982). Editors, S. S.
Chern and W. T. Wu.

\bibitem{5}
K. S. Chou and X. Wang, {\em A variational theory of the Hessian
equation}, { Comm. Pure Applied Math.}, {\bf 54} (2001),
1029-1064.

\bibitem{6}
B. Guan, {\em The Dirichlet problem for a class of fully nonlinear
elliptic equations}, {Comm. Partial Diffrerntial Equations}, {\bf
19} (1994), 399-416.

\bibitem{7}
B. Guan and H. Y. Jian, {\em The Monge-Amp\`ere equation with
infinite boundary value}, Pacific J. Math., {\bf 216} (2004),
77-94.

\bibitem{8}
B. Guan and P. Guan, {\em Convex hypersurfaces of prescribed
curvature}, {Annal of Math.}, {\bf 156} (2002), 655-674.

\bibitem{9}
P. Guan, {\em Geometric fully nonlinear equations}, {Lectures at
CMS, Zhejiang Univ.}, 2004.

\bibitem{10}
P. Guan and X. N. Ma, {\em Christoffel-Minkowski problem},
{Inventions Math.}, {\bf 151} (2003), 553-577.

\bibitem{11}
J. B. Keller, {\em On solutions of $\Delta u = f (u)$}, {Comm.
Pure Applied Math.}, {\bf 10} (1957), 503-510.

\bibitem{12}
P. L. Lions, {\em  Sur les \'equations de Monge-Amp\`ere I},
  Manuscripta Math., {\bf 41} (1983), 1-43;
{\em II}, Arch. Rational Mech. Anal., {\bf 89} (1985), 93-122.

\bibitem{13}
J. Matero, {\em The Bieberbach-Rsdemacher problem for the
Monge-Amp\`ere operator}, Manuscripta Math., {\bf 91} (1996),
379-391.

\bibitem{14}
R. Osserman, {\em On the inequality $\Delta u \geq f (u)$},
Pacific J. Math., {\bf 7} (1957), 1641-1647.

\bibitem{15}
P. Salani, {\em Boundary blow-up problems for Hessian equations} ,
Manuscripta Math., {\bf 96} (1998), 281-294.

\bibitem{16}
N. S. Trudinger, {\em On the Dirichlet problem for Hessian
equations}, {Acta Math.}, {\bf 175} (1995), 151-164.

\bibitem{17}
N. S. Trudinger, {\em Weak solutions of Hessian equations}, Comm.
Partial Diff. Eqns., {\bf 22} (1997), 1251-1651.

\bibitem{18}
J. I. E. Urbas, {\em On the existence of nonclassial solutions for
two classes of fully nonlinear elliptic equations}, { India Univ.
Math. J.}, {\bf 39} (1990), 355-382.


\end{thebibliography}
\end{document}